# LARGE DEVIATIONS FOR RANDOM POWER MOMENT PROBLEM

By Fabrice Gamboa and Li-Vang Lozada-Chang

*Université Paul Sabatier and Universidad de la Habana*

We consider the set $M_n$ of all $n$-truncated power moment sequences of probability measures on $[0,1]$. We endow this set with the uniform probability. Picking randomly a point in $M_n$, we show that the upper canonical measure associated with this point satisfies a large deviation principle. Moderate deviation are also studied completing earlier results on asymptotic normality given by Chang, Kemperman and Studden [*Ann. Probab.* **21** (1993) 1295–1309]. Surprisingly, our large deviations results allow us to compute explicitly the $(n+1)$th moment range size of the set of all probability measures having the same $n$ first moments. The main tool to obtain these results is the representation of $M_n$ on canonical moments [see the book of Dette and Studden].

**1. Introduction.** In this work we will study the asymptotic behavior in large deviations of random power moment problem. Let $\mathbb{P}([0,1])$ denote the set of all probability measures (p.m.s) on the interval $[0,1]$. In the whole paper this set will be endowed with the weak topology [see Billingsley (1999)]. For any $\mu$ in $\mathbb{P}([0,1])$ the $k$th (power) moment of $\mu$ will be denoted by $c_k(\mu)$:

$$c_k(\mu) = \int_{[0,1]} x^k \, d\mu(x).$$

In this paper we focus on some asymptotic properties of the finite moment space $M_n$:

$$M_n = \{c^{(n)}(\mu) = (c_1(\mu), \ldots, c_n(\mu))^T : \mu \in \mathbb{P}([0,1])\}, \qquad n \in \mathbb{N}^*.$$

$M_n$ is the closed convex hull of the curve

$$\{(x, x^2, \ldots, x^n)^T : x \in [0,1]\}$$

[see Karlin and Studden (1966) and Kreĭn and Nudelman (1977)]. As $M_n$ is a compact subset of $\mathbb{R}^n$ having nonvoid interior, we may define the uniform probability $\mathbb{P}_n$ on $M_n$. $M_n$ is a very "small" set. Indeed, its volume has order $2^{-n^2}$ for large $n$. Moreover, this set "concentrates" (in some sense) on a single point. More precisely, Chang, Kemperman and Studden (1993) have shown that for any fixed $k \in \mathbb{N}^*$ and under the probability $\mathbb{P}_n$, the $k$ first components of $c^{(n)} \in M_n$ converge in probability to the $k$ first moments of the arcsine law. That is,

$$(1) \qquad \lim_{n \to \infty} \mathbb{P}_n \left( \left\| \begin{pmatrix} c_1^{(n)} \\ \vdots \\ c_k^{(n)} \end{pmatrix} - \begin{pmatrix} \bar{c}_1 \\ \vdots \\ \bar{c}_k \end{pmatrix} \right\| \geq \xi \right) = 0, \qquad \xi > 0,$$

where $\bar{c}_j = c_j(\nu)$, $j \in \mathbb{N}^*$, and

$$(2) \qquad \nu(dx) = \frac{dx}{\pi \sqrt{x(1-x)}}.$$

Hence, $M_n$ is concentrated around the $n$th first moments of the arcsine law. Moreover, Chang, Kemperman and Studden (1993) have studied the fluctuation limit law in (1). They have shown that the limit distribution of the fluctuations is Gaussian:

$$(3) \qquad \sqrt{n}(Z_n^{(k)} - \bar{c}^{(k)}) \xrightarrow[n \to \infty]{\mathcal{L}-\mathbb{P}_n} \mathcal{N}_k(0, \Sigma_k),$$

where $Z_n^{(k)}$ denotes the random vector built with the $k$ first coordinates of $Z_n = c^{(n)}$ (drawn randomly with distribution $\mathbb{P}_n$). $\bar{c}^{(k)} = (\bar{c}_1, \bar{c}_2, \ldots, \bar{c}_k)^T$ is the vector of the $k$ first moment of $\nu$. The covariance matrix $\Sigma_k$ will be described in Section 2.5. The main result of this paper is a functional large deviations principle for the sequence $(Z_n^{(k)})$. This means that $Z_n^{(k)}$ "concentrates" exponentially fast. More precisely, for any $c^{(n)} \in M_n$, there exists a unique measure $\sigma_n^+(c^{(n)})$ whose $n$ first moments are $c^{(n)}$ and maximizing the $(n+1)$th moment [see Section 2.3, Karlin and Studden (1966), Kreĭn and Nudelman (1977) and Dette and Studden (1997)]. In Section 2.4 we show that $\sigma_n^+(Z_n)$ satisfies a large deviations principle (LDP). Roughly speaking, this means that $\sigma_n^+(Z_n)$ concentrates exponentially fast on $\nu$. In other words, for large $n$,

$$(4) \qquad \mathbb{P}_n(\sigma_n^+(Z_n) \in A) \approx \exp\left(-n \inf_{\mu \in A} I(\mu)\right),$$

where $A$ is a Borel-measurable set of $\mathbb{P}([0,1])$ and $I$ is the rate function of the LDP (see Sections 2.1 and 2.4 for the more precise statement). Although we prove this main result using a limit projective approach [see Dembo and Zeitouni (1998)] we manage to compute precisely the rate function $I$. This is performed using Szegö asymptotic theory on orthogonal polynomials [Grenander and Szegö (1958)]. Surprisingly, $I$ is the reversed

Kullback information (or cross entropy) with respect to $\nu$. Moreover, this computation gives quantitative evaluations of the $(n+1)$th moment range size of the set of all probability measures having the same prescribed $n$ first moments in term of the reversed Kullback information. To our knowledge these results on power moment problem are new, they are developed in Section 3.

In Section 2.5, we give large and moderate deviation principles for the random vector $Z_n^{(k)}$. Hence, we study exponential rates of convergence, between (3) and (4), for this vector. All the proofs are postponed to Section 4. This last section begin with a section on canonical moments [Dette and Studden (1997)] which are the main tool to prove our results.

## 2. Main results.

2.1. *Large deviations.* Let us first recall what is a LDP [see, e.g., Dembo and Zeitouni (1998)]. Let $(u_n)$ be a decreasing positive sequence of real numbers.

DEFINITION 2.1. We say that a sequence $(R_n)$ of probability measures on a measurable Hausdorff space $(U, \mathcal{B}(U))$ satisfies a LDP with rate function $I$ and speed $(u_n)$ if:

(i) $I$ is lower semicontinuous (l.s.c.), with values in $\mathbb{R}^+ \cup \{+\infty\}$.
(ii) For any measurable set $A$ of $U$,

$$-I(\operatorname{int} A) \leq \liminf_{n \to \infty} u_n \log R_n(A) \leq \limsup_{n \to \infty} u_n \log R_n(A) \leq -I(\operatorname{clo} A),$$

where $I(A) = \inf_{\xi \in A} I(\xi)$ and $\operatorname{int} A$ (resp. $\operatorname{clo} A$) is the interior (resp. the closure) of $A$.

We say that the rate function $I$ is good if its level set $\{x \in U : I(x) \leq a\}$ is compact for any $a \geq 0$. More generally, a sequence of $U$-valued random variables is said to satisfy a LDP if their distributions satisfy a LDP.

To be self-contained let us recall some facts and tools on large deviations which will be useful in the paper [we refer to Dembo and Zeitouni (1998) for more on large deviations]:

- *Contraction principle.* Assume that $(R_n)$ satisfies a LDP on $(U, \mathcal{B}(U))$ with good rate function $I$ and speed $(u_n)$. Let $T$ be a continuous mapping from $U$ to another space $V$. Then, $(R_n \circ T^{-1})$ satisfies a LDP on $(V, \mathcal{B}(V))$ with good rate function

$$I'(y) = \inf_{x : T(x) = y} I(x), \qquad y \in V,$$

and speed $(u_n)$.

- *Exponential approximation.* Assume that $U$ is a metric space and let $d$ denote the distance on $U$. Let $(X_n)$ be a $U$-valued random sequence satisfying a LDP with good rate function $I$ and speed $u_n$. Let $(Y_n)$ be another $U$-valued random sequence. If for any $\xi > 0$,

$$\limsup_{n \to \infty} u_n \log \mathbb{P}(d(X_n, Y_n) > \xi) = -\infty,$$

then $(Y_n)$ shares the same LDP as $(X_n)$.

In the sequel, talking about large deviations, if we omit the sequence $(u_n)$ it means that $u_n = \frac{1}{n}$ ($n \in \mathbb{N}^*$). Following Section 3.7 of Dembo and Zeitouni (1998), when the sequence $(u_n)$ satisfies $n^{-1} = o(u_n)$ we will say moderate deviations rather than large deviations with speed $(u_n)$.

2.2. *Kullback and reversed Kullback information.* Let $(U, \mathcal{B}(U))$ be a measurable space and $P$ and $Q$ be p.m.s. on $(U, \mathcal{B}(U))$. Recall that the *Kullback information* or cross entropy of $P$ with respect to $Q$ is defined by

$$K(P, Q) = \begin{cases} \int_U \log \frac{dP}{dQ} dP, & \text{if } P \ll Q \text{ and } \log \frac{dP}{dQ} \in L^1(P), \\ +\infty, & \text{otherwise.} \end{cases}$$

Properties of $K$ as a function of $P$ may be found in Bretagnolle (1979). $K$ is the rate function for Sanov large deviations theorem [see Dembo and Zeitouni (1998)]. The rate function involved in this paper is the reversed Kullback information with respect to $\nu$, that is,

$$I(\mu) = K(\nu, \mu), \qquad \mu \in \mathbb{P}([0,1]),$$

where $\nu$ is the arcsine law (2). Observe that $I$ is l.s.c. [Theorem 2.7 of Borwein and Lewis (1993)]. Moreover, it is obviously a good rate function [$\mathbb{P}([0,1])$ is a compact set].

The following property may be found in Theorem 2.1 of Gamboa and Gassiat (1997) and in Section 3 of Borwein and Lewis (1993):

PROPOSITION 2.2. (i) *For* $\mu \in \mathbb{P}([0,1])$,

$$I(\mu) = \sup_{f \in C[0,1]} \left( \int_0^1 f(x) \, d\mu(x) + \int_0^1 \ln(1 - f(x)) \, d\nu(x) \right),$$

*where $C[0,1]$ is the set of all continuous functions on $[0,1]$. (In the whole of the paper we take the convention $\ln \tau = -\infty$ whenever $\tau \leq 0$.)*

(ii) *Let $\Phi = (\Phi_1, \ldots, \Phi_k)^T \in (C[0,1])^k$ and for $c \in \mathbb{R}^k$,*

$$\mathbb{S}_\Phi(c) = \left\{ \mu \in \mathbb{P}([0,1]) : \int_0^1 \Phi_j(x) \, d\mu(x) = c_j, j = 1, \ldots, k \right\}.$$

*Then,*

$$\inf_{\mu \in \mathbb{S}_\Phi(c)} I(\mu) = \sup_{(\lambda_0, \lambda) \in \mathbb{R}^{k+1}} \left\{ \lambda_0 + \langle \lambda, c \rangle + \int_0^1 \ln(1 - \lambda_0 - \langle \lambda, \Phi(x) \rangle) \, d\nu(x) \right\},$$

*where $\langle \cdot, \cdot \rangle$ denotes the usual scalar product on $\mathbb{R}^k$.*

We will come back on the second point in Section 2.5. This property will be helpful to show the results of Section 3.

For $x \in [0,1]$ and $k \in \mathbb{N}^*$, let $\phi_k(x) = (x, x^2, \ldots, x^k)^T$. For $c \in \mathbb{R}^k$, we set $\mathbb{S}_k(c) = \mathbb{S}_{\phi_k}(c)$. In the next section we briefly recall some known facts on $\mathbb{S}_k(c)$. In Section 3 we give new results on $\mathbb{S}_k(c)$.

2.3. *Power moment problem.* Let $M$ denote the set of all infinite moment sequences, that is,

$$M = \{c^{(\infty)}(\mu) = (c_i(\mu))_{i \in \mathbb{N}} : \mu \in \mathbb{P}([0,1])\}.$$

For $n \in \mathbb{N}^*$, let $\Pi_n : M \to M_n$ denote the natural projection map. So, we have $M_n = \Pi_n(M)$.

Let recall some useful facts on the power moments problems [see Karlin and Studden (1966) and Kreĭn and Nudelman (1977) for an exhaustive overview on this problem]. Because $[0,1]$ is a compact interval, the relation between the elements of $M$ and $\mathbb{P}([0,1])$ is bijective. In general, for $c^{(n)} \in M_n$, there exists an infinite number of probability measures such that $c^{(n)} = c^{(n)}(\mu)$. More precisely, for any $n \in \mathbb{N}^*$ and $c^{(n)} \in \mathbb{R}^n$:

(i) $\#\mathbb{S}_n(c^{(n)}) = +\infty \Leftrightarrow c^{(n)} \in \operatorname{int} M_n$,
(ii) $\#\mathbb{S}_n(c^{(n)}) = 1 \Leftrightarrow c^{(n)} \in \partial M_n$ (the boundary of $M_n$),
(iii) $\#\mathbb{S}_n(c^{(n)}) = 0 \Leftrightarrow c^{(n)} \notin M_n$,

where $\#A$ is the number of elements lying in $A$. These results come from the fact that $(\phi_n)$ is a Tchebycheff system [see Karlin and Studden (1966) and Kreĭn and Nudelman (1977)]. The elements of $\partial M_n$ satisfy an extremal property. For $c^{(k)} \in M_k$, let $c^+, c^- : M_k \to \mathbb{R}$ be defined by

$$c^+(c^{(k)}) = \max\{c \in \mathbb{R} : (c_1, c_2, \ldots, c_k, c)^T \in M_{k+1}\},$$

$$c^-(c^{(k)}) = \min\{c \in \mathbb{R} : (c_1, c_2, \ldots, c_k, c)^T \in M_{k+1}\}.$$

Then, $c^{(k+1)} = ((c^{(k)})^T, c_{k+1})^T \in \partial M_{k+1}$ if, and only if, $c_{k+1} = c^+(c^{(k)})$ [or $c^-(c^{(k)})$].

Let $d^{(n)} \in M_n$. We will denote by $\sigma_n^+(d^{(n)})$ the measure $\mu$ such that $c^{(n)}(\mu) = d^{(n)}$ and $c_{n+1}(\mu) = c^+(d^{(n)})$. This measure is the so-called upper canonical representation of the finite moments sequence $d^{(n)}$ [see Karlin and Studden (1966) and Kreĭn and Nudelman (1977)]. In the next section we will study large deviations properties of $\sigma_n^+(Z_n)$.

2.4. *Large deviations for $\sigma_n^+(Z_n)$.* Recall that $\mathbb{P}_n$ denotes the normalized Lebesgue probability measure on $M_n$ and $Z_n$ denotes a random vector having distribution $\mathbb{P}_n$. For $k \leq n$, let $\Pi_k^n : M_n \to M_k$ be the projection map taking the $k$ first coordinates of an element of $M_n$. Let $Z_n^{(k)}$ denote the random vector $\Pi_k^n(Z_n)$ and $\mu_n$ be the random measure $\sigma_n^+(Z_n)$. Let $Q_n$ be the law of $\mu_n$. Our main result is the large deviations principle for $(\mu_n)$.

THEOREM 2.3. *$(\mu_n)$ satisfies a LDP with convex good rate function $I$.*

COROLLARY 2.4. *$(\mu_n)$ converges in probability to $\nu$.*

We can also study the LDP associated to the so called tilted measures. Changing a little the measures $Q_n$ on $\mathbb{P}([0,1])$, we obtain the new limit law for the sequence of random probability measures $\mu_n$. Namely, consider on $\mathbb{P}([0,1])$ the new probability measure $\widetilde{Q}_n$ defined by

$$\widetilde{Q}_n(B) = \frac{\mathbb{E}_{Q_n}(\exp(nF(\mu_n))\mathbb{1}_B(\mu_n))}{\mathbb{E}_{Q_n} \exp(nF(\mu_n))}, \tag{5}$$

where $B$ runs over the Borelian sets of $\mathbb{P}([0,1])$ and $F : \mathbb{P}([0,1]) \to \mathbb{R}$ is the continuous functional defined by

$$F(\mu) = \int_{[0,1]} f_0(x) \, d\mu(x),$$

where $f_0 \in C([0,1])$.

We have the following results.

THEOREM 2.5. *$(\widetilde{Q}_n)_{n \in \mathbb{N}}$ satisfies a LDP on $\mathbb{P}([0,1])$ with good rate function*

$$I_F(\mu) = I(\mu) - F(\mu) + K_F,$$

*where $K_F := \sup_{\mu' \in \mathbb{P}([0,1])} \{F(\mu') - I(\mu')\}$.*

REMARK 2.1. Let $\tilde{\mu}_n$ denote a random measure of $\mathbb{P}([0,1])$ having distribution $\widetilde{Q}_n$. The existence of an unique minimum point for $I_F$ implies the convergence (in probability) of $\tilde{\mu}_n$ toward this minimum point. Under certain conditions over $f_0$ we can characterize the minimum points of $I_F$. Let $\underline{\mu} \in \mathbb{P}([0,1])$ be a minimum point of $I_F$. Let $\underline{\mu} = \underline{\mu}_a + \underline{\mu}_s$ be the Lebesgue decomposition of $\underline{\mu}$ with respect to $\nu$ ($\underline{\mu}_a \ll \nu$). Then, from Theorem 3.5 of Borwein and Lewis (1993) there exists $\lambda^* \in \mathbb{R}$ with:

(i) $g_{\lambda^*} = \frac{d\underline{\mu}_a}{d\nu} = \frac{1}{\lambda^* - f_0}$ a.s.,
(ii) $\operatorname{supp} \underline{\mu}_s \subset \{x \in [0,1] : f_0(x) = \lambda^*\}$.

We now give two particular cases where the sequence $(\tilde{\mu}_n)$ has a limit. Let $\bar{\lambda} := \max_{x \in [0,1]} f_0(x)$ and $\chi_0 := \int_{[0,1]} (\bar{\lambda} - f_0)^{-1} d\nu$.

(i) Assume $\chi_0 \geq 1$, then there exist $\lambda^* \geq \bar{\lambda}$ such that
$$\int_{[0,1]} g_{\lambda^*}(x) \, d\nu(x) = \int_{[0,1]} \frac{d\nu(x)}{\lambda^* - f_0(x)} = 1$$
and $(\tilde{\mu}_n)$ converges in probability to $g_{\lambda^*}\nu$.

(ii) If $\chi_0 < 1$ and $\{x \in [0,1] : f_0 = \bar{\lambda}\}$ reduces to the singleton $\{x_0\}$, then $(\tilde{\mu}_n)$ converges in probability to
$$g_{\bar{\lambda}}\nu + (1 - \chi_0)\delta_{x_0}.$$

2.5. *Large and moderate deviations for finite moments sequence.* Although we may obtain a LDP for $(Z_n^{(k)})$ using Theorem 2.3 and the contraction principle, we will show the following theorem.

THEOREM 2.6. $(Z_n^{(k)})$ *satisfies a LDP with convex good rate function*

(6) $$I_k(c^{(k)}) = \begin{cases} -\ln(c^+(c^{(k)}) - c^-(c^{(k)})) - k\ln 4, & \text{if } c^{(k)} \in \text{int } M_k, \\ +\infty, & \text{otherwise.} \end{cases}$$

Theorem 2.3 will be shown by first proving Theorem 2.6 and a projective limit argument [see Section 4.6 in Dembo and Zeitouni (1998)].

REMARK 2.2. $I_k$ achieves the value 0 only at $\bar{c}^{(k)}$. Hence, using Borel Cantelli lemma, that yields to the almost sure convergence of $Z_n^{(k)}$ to $\bar{c}^{(k)}$.

We now turn on moderate deviation properties of $(Z_n^{(k)})$. For $i, j \in \mathbb{N}^*$, define

(7) $$a_{ij} = \begin{cases} 2^{-2i+1}\binom{2i}{i-j}, & \text{if } 1 \leq j \leq i, \\ 0, & \text{if } j > i. \end{cases}$$

For $k \in \mathbb{N}^*$, let $A_k = (a_{ij})_{i,j=1}^k$ and set
$$\Sigma_k = \tfrac{1}{2} A_k A_k^T$$
and
$$J_k(\bar{x}) = \tfrac{1}{2} \bar{x}^T \Sigma_k^{-1} \bar{x}, \qquad \bar{x} \in \mathbb{R}^k.$$

THEOREM 2.7. *Let $(u_n)$ be a sequence decreasing to 0 such that $n^{-1} = o(u_n)$ and let $\widetilde{Z}_n^{(k)} := \sqrt{nu_n}(Z_n^{(k)} - \bar{c}^{(k)})$. Then $(\widetilde{Z}_n^{(k)})$ satisfies a moderate deviations principle with good rate function $J_k$.*

**3. On the $(n+1)$th moment range size of the probability measures having the same $n$ first moments.** In this section, for $c^{(k)} \in M_k$ ($k \in \mathbb{N}^*$), we set

$$r_{k+1}(c^{(k)}) = \sup_{\mu_1, \mu_2 \in \mathbb{S}_k(c^{(k)})} \left\{ \int_0^1 x^{k+1} \, d\mu_1(x) - \int_0^1 x^{k+1} \, d\mu_2(x) \right\}.$$

Using the large deviations properties of the previous section and the contraction principle, we obviously obtain the following:

THEOREM 3.1. *Let $c^{(k)} \in M_k$. Then*:

(i)

(8) $$r_{k+1}(c^{(k)}) = \exp\left(-\inf_{\mu \in \mathbb{S}_k(c^{(k)})} I(\mu) - k \ln 4\right).$$

(ii) *Let $\mu \in \mathbb{P}([0,1])$ and $c^{(k)} = c^{(k)}(\mu)$, then*

$$\lim_{k \to \infty} \frac{1}{k} \ln[4^k r_{k+1}(c^{(k)})] = -I(\mu).$$

REMARK 3.1. (i) Using Proposition 2.2, (8) may be expressed as the supremum of a concave function. Indeed, we have

$$r_{k+1}(c^{(k)}) = \exp\left(-k \ln 4 - \sup_{(\lambda_0, \lambda) \in \mathbb{R}^{n+1}} H_k(\lambda_0, \lambda, c^{(k)})\right),$$

where for $(\lambda_0, \lambda) \in \mathbb{R}^{k+1}$,

$$H_k(\lambda_0, \lambda, c^{(k)}) = \lambda_0 + \langle \lambda, c^{(k)} \rangle + \int_0^1 \ln(1 - \lambda_0 - \langle \lambda, \phi_k(x) \rangle) \, d\nu(x).$$

So for all $(\lambda_0, \lambda) \in \mathbb{R}^{k+1}$,

$$r_{k+1}(c^{(k)}) \leq \exp[-k \ln 4 - H_k(\lambda_0, \lambda, c^{(k)})].$$

This last inequality is helpful to study the superresolution rate in the power moment problem. This will be done in a forthcoming paper of Gamboa and Lozada. Superresolution occurs for $c^{(k)} \in \partial M_k$. In this case as we saw in Section 2.3, $\mathbb{S}_k(c^{(k)})$ reduces to a single p.m. Superresolution rate is the concentration rate of the set $\mathbb{S}_k(c^{(k)} + \eta)$ when $\eta \in \mathbb{R}^k$ is a small perturbation [see Gamboa and Gassiat (1996) and Doukhan and Gamboa (1996) for more on this problem].

(ii) Let $P$ be a polynomial having degree $k \in \mathbb{N}^*$. Assume that $P$ is positive on $[0,1]$ and satisfies the normalizing condition

$$\int_0^1 \frac{d\nu(x)}{P(x)} = 1.$$

Set, for $j \in \mathbb{N}^*$, $d_j = c_j(\frac{\nu}{P})$ and, as usual, $d^{(j)} = (d_i)_{i=1,\ldots,j}$. Then, using optimization results developed in Theorem 2.1 of Gamboa and Gassiat (1997) or Theorem 4.1 of Borwein and Lewis (1993) we have, for $j \geq k$,

$$I\left(\frac{\nu}{P}\right) = \inf_{\mu \in \mathbb{S}_j(d^{(j)})} I(\mu)$$

(see also Section 2.5 for related results). Therefore, in this frame Theorem 3.1 gives

$$r_{j+1}(d^{(j)}) = \exp\left(-j\ln 4 - \int_0^1 \ln P(x)\,d\nu(x)\right), \qquad j \geq k.$$

(iii) Obviously, using the contraction principle and Theorem 2.6, we obtain that the sequence $(r_{k+1}(Z_n^{(k)}))$ satisfies a LDP with good rate function,

$$R(r) = \begin{cases} -\ln(4^k r), & \text{if } r \in\, ]0, 4^{-k}[, \\ +\infty, & \text{otherwise.} \end{cases}$$

## 4. Proofs.

4.1. *Canonical moments.* This section is devoted to canonical moments which are the basic tool for our results [as those obtained in Chang, Kemperman and Studden (1993)]. We will present few properties of canonical moments. An exhaustive study of canonical moments and their applications can be found in the very nice book of Dette and Studden (1997). Let $p_i$ denote the $i$th canonical moments defined in $\text{int}\, M_k$, $k \in \mathbb{N}^*$, by

$$p_i(c^{(k)}) = p_i(c^{(i)})$$
$$= \frac{c_i - c^-(c^{(i-1)})}{c^+(c^{(i-1)}) - c^-(c^{(i-1)})}, \qquad 1 \leq i \leq k,$$

where $c^+$ and $c^-$ have been defined in Section 2.3. For $k \in \mathbb{N}^*$, we will denote by $\mathbf{p}_k$ the map from $\text{int}\, M_k$ to $]0,1[^k$,

$$c^{(k)} \mapsto p^{(k)}(c^{(k)}) = (p_1(c_1), \ldots, p_2(c^{(2)}), \ldots, p_k(c^{(k)}))^T.$$

The map $\mathbf{p}_k$ is a continuous one-to-one onto correspondence between the interior of the $k$th moment space and $]0,1[^k$ [see, Skibinsky (1968) and Dette and Studden (1997)]. Obviously, the sequence $(\mathbf{p}_k)$ induces a one-to-one correspondence $\mathbf{p}_\infty$ between $\text{int}\, M$ and $]0,1[^{\mathbb{N}}$. A very interesting property of the canonical moments is that they are invariants under linear transformations of the interval. Therefore, as the reader will easily check, all the results obtained on $[0,1]$ are also valid for any other bounded closed interval of the real line. The canonical moments of the arcsine law $\nu$ are all equal to $\frac{1}{2}$ [see Skibinsky (1969)]. We will denote them by $\bar{p}_j$.

In Lemma 1.4 of Chang, Kemperman and Studden (1993) the first Taylor expansion of $\mathbf{p}_k$ around $\bar{p}^{(k)}$ is given,

$$(9) \quad c_m = \bar{c}_m + 2\sum_{j=1}^m a_{mj}(p_j - \tfrac{1}{2}) + \mathbf{O}\left(\sum_{j=1}^m |p_j - \tfrac{1}{2}|^2\right),$$

where, for $m, j \in \mathbb{N}^*$, $c^{(m)} = (c_1, \ldots, c_m)^T \in \text{int}\, M_m$, $p^{(m)}(c^{(m)}) = (p_1, \ldots, p_m) \in ]0,1[^m$ and $a_{mj}$ has been defined in Section 2.5.

We recall that the beta distribution $(\beta(a,b))$ of parameters $a, b > 0$, has density with respect to the Lebesgue measure on $[0,1]$,

$$[B(a,b)]^{-1} x^{a-1}(1-x)^{b-1},$$

where $B(a,b) = \int_{[0,1]} x^{a-1}(1-x)^{b-1}\, dx$.

We have the following lemma.

LEMMA 4.1 [Theorem 1.3 of Chang, Kemperman and Studden (1993)]. *If $M_n$ is endowed with $\mathbb{P}_n$, then the random vector $p^{(n)} = (p_1, p_2, \ldots, p_n)^T$ satisfies:*

  (i) $(p_i)_{i=1}^n$ *are independent random variables,*
  (ii) $p_i \sim \beta(n-i+1, n-i+1), i = 1, \ldots, n$.

Now, we display an equation linking canonical and ordinary power moments given in Skibinsky (1967), for $c^{(n)} \in \text{int}\, M_n$,

$$(10) \quad r_{n+1}(c^{(n)}) = c^+(c^{(n)}) - c^-(c^{(n)}) = \prod_{i=1}^n p_i(c^{(n)})(1 - p_i(c^{(n)})).$$

In (15) we will express $r_{n+1}(c^{(n)})$ as a function of Hankel determinants.

4.2. *Proof of Theorem* 2.6.

LEMMA 4.2 [Exercise 4.2.7 of Dembo and Zeitouni (1998)]. *Let $(X_{n,1})$ [resp. $(X_{n,2})$] be a sequence of random variables taking its values in a regular space $\mathcal{X}_1$ (resp. $\mathcal{X}_2$) satisfying a LDP with good rate function $I_1$ (resp. $I_2$). If $X_{n,1}$ and $X_{n,2}$ are independent, then $X_n = (X_{n,1}, X_{n,2})$ satisfies a LDP in $\mathcal{X}_1 \times \mathcal{X}_2$ with good rate function $I(x_1, x_2) = I_1(x_1) + I_2(x_2)$.*

LEMMA 4.3. *The sequence of p.m.s $(\beta(n,n))_{n \in \mathbb{N}}$ satisfies on $[0,1]$ a LDP with good rate function*

$$\hat{I}(x) = -\ln(x - x^2) - \ln 4.$$

PROOF. Let $X_n$ be a random variable having $\beta(n,n)$ distribution, then we may write

$$X_n \stackrel{\mathcal{L}}{=} \frac{\sum_{i=1}^n Y_i}{\sum_{i=1}^{2n} Y_i},$$

where the random variables $(Y_i)$ are independent with standard exponential distribution [see, e.g., Bartoli and Del Moral (2001), page 71]. Consider the bidimensional random vector

$$V_n = \left(\frac{1}{n}\sum_{i=1}^n Y_i, \frac{1}{n}\sum_{i=n+1}^{2n} Y_i\right)^T.$$

From Lemma 4.2 and Cramér's theorem [Theorem 2.2.3 in Dembo and Zeitouni (1998)], we obtain a LDP for $(V_n)$ with good rate function

$$I_V(x,y) = x + y - 2 + \log xy, \qquad x > 0, y > 0.$$

Hence, the LDP for $(X_n)$ follows from the one for $(V_n)$ and the contraction principle with the continuous map $(x,y) \mapsto x/(x+y)$. □

From Lemmas 4.3 and 4.1, we have that, for $j \in \mathbb{N}^*$, $(p_j(Z_n))_{n \in \mathbb{N}}$ ($j \in \mathbb{N}$, fixed) satisfies a LDP with good rate function $\hat{I}$. Thanks to Lemma 4.2, for $k \in \mathbb{N}$, we also have a LDP for $(p^{(k)}(Z_n))_{n \in \mathbb{N}}$ with good rate function

$$\hat{I}_k(p^{(k)}) = -\sum_{i=1}^k \ln(p_i - p_i^2) - k \ln 4,$$

where $p^{(k)} = (p_1, p_2, \ldots, p_k)^T \in ]0,1[^k$.

We have $Z_n^{(k)} = \mathbf{p}_k^{-1}(p^{(k)}(Z_n))$. The contraction principle yields the LDP for $(Z_n^{(k)})_{n \in \mathbb{N}}$ with good rate function,

$$I_k(c^{(k)}) = \hat{I}_k(p^{(k)}(c^{(k)})).$$

Using (10), we may write $I_k$ as in (6).

It is obvious that the function $\hat{I}_k$ achieves its minimum value 0 at $\bar{p}^{(k)}$ (and only at this point). Indeed, remember that $\bar{p}_j = \frac{1}{2}$ for $1 \leq j \leq k$. Consequently, the function $I_k$ achieves its minimum value at $\bar{c}^{(k)}$.

4.3. *Proof of Theorem 2.7.*

LEMMA 4.4. *Let $(u_n)$ be a sequence decreasing to 0 and $n^{-1} = o(u_n)$. Let $Y_n$ be a random variable having $\beta(n-l, n-l)$ distribution $(n \in \mathbb{N}^*, l \in \mathbb{N})$. Set*

$$X_n := \sqrt{nu_n}(Y_n - \tfrac{1}{2}),$$

*then $(X_n)$ satisfies a LDP with good rate function $J_1(x) = 4x^2$ and speed $(u_n)$.*

PROOF. It is well known that $B_n = ((n-1)!)^2/(2n-1)!$. Using Stirling's formula for $n!$, we have

$$(11) \qquad B_n = 4^{-n} \left(\frac{\pi}{n}\right)^{1/2} (2 + \xi_n)$$

with $\xi_n \to 0$ as $n \to \infty$. Therefore, $\lim_n \frac{1}{n} \ln B_n = -\ln 4$.

Let $\mathcal{A} := \{]a, b[ : a < b\}$. $\mathcal{A}$ is a base of the usual topology on $\mathbb{R}$. Let $]a, b[ \in \mathcal{A}$. For $n$ large enough,

$$\mathbb{P}(X_n \in ]a, b[) = \mathbb{P}\left(Y_n \in \left]\frac{1}{2} + \frac{a}{\sqrt{nu_n}}, \frac{1}{2} + \frac{b}{\sqrt{nu_n}}\right[\right)$$

$$= B_{n-l-1}^{-1} \int_{1/2 + a/\sqrt{nu_n}}^{1/2 + b/\sqrt{nu_n}} x^{n-l-1}(1-x)^{n-l-1} \, dx$$

$$\leq B_{n-l-1}^{-1} m^+\left(\frac{1}{2} + \frac{a}{\sqrt{nu_n}}, \frac{1}{2} + \frac{b}{\sqrt{nu_n}}\right)^{n-l-1} \frac{b-a}{\sqrt{nu_n}},$$

where $m^+(a, b) := \sup\{t - t^2 : t \in ]a, b[\}$ and $m^-(a, b) := \inf\{t - t^2 : t \in ]a, b[\}$. The function $t \mapsto t^2 - t$ is strictly convex in $]0, 1[$ with minimum value at $\frac{1}{2}$, so

$$m^+\left(\frac{1}{2} + \frac{a}{\sqrt{nu_n}}, \frac{1}{2} + \frac{b}{\sqrt{nu_n}}\right) = \begin{cases} \frac{1}{4} - \frac{a^2}{nu_n}, & \text{if } a > 0, \\ \frac{1}{4} - \frac{b^2}{nu_n}, & \text{if } b < 0, \\ 0, & \text{if } a < 0 < b. \end{cases}$$

Using (11),

$$\limsup_n u_n \ln \mathbb{P}(X_n \in ]a, b[) \leq \begin{cases} -4a^2, & \text{if } a > 0, \\ -4b^2, & \text{if } b < 0, \\ 0, & \text{if } a < 0 < b. \end{cases}$$

This leads to

$$(12) \qquad \sup_{\{A \in \mathcal{A} : x \in A\}} \left\{ -\limsup_n u_n \ln \mathbb{P}(X_n \in A) \right\} \geq 4x^2.$$

Now, analogously, we have

$$\mathbb{P}(X_n \in A) \geq m^-\left(\frac{1}{2} + \frac{a}{\sqrt{nu_n}}, \frac{1}{2} + \frac{b}{\sqrt{nu_n}}\right)^{n-l-1} \frac{b-a}{\sqrt{nu_n}} \frac{1}{B_{n-l-1}}.$$

But,

$$m^-\left(\frac{1}{2}+\frac{a}{\sqrt{nu_n}},\frac{1}{2}+\frac{b}{\sqrt{nu_n}}\right) = \begin{cases} \frac{1}{4}-\frac{b^2}{nu_n}, & \text{if } a>0, \\ \frac{1}{4}-\frac{a^2}{nu_n}, & \text{if } b<0, \\ \frac{1}{4}-\frac{\max\{a^2,b^2\}}{nu_n}, & \text{if } a<0<b. \end{cases}$$

So, using (11), we obtain

$$\limsup_n u_n \ln \mathbb{P}(X_n \in ]a,b[) \leq \begin{cases} -4b^2, & \text{if } a>0, \\ -4a^2, & \text{if } b<0, \\ -4\max\{a^2,b^2\}, & \text{if } a<0<b. \end{cases}$$

Consequently,

(13) $$\sup_{\{A\in\mathcal{A}\,:\,x\in A\}}\left\{-\liminf_n u_n \ln \mathbb{P}(X_n \in A)\right\} \leq 4x^2.$$

Using Theorem 4.1.11 of Dembo and Zeitouni (1998), (12) and (13), we get a weak LDP for $X_n$ with rate function $J_1$. The full LDP is then a consequence of the exponential tightness of the laws of $Y_n$. Indeed, for $K>0$ and $n$ large enough, we have

$$\mathbb{P}(X_n \in [-K,K]^c) = 2\mathbb{P}\left(Y_n \in \left]\frac{1}{2}+\frac{K}{\sqrt{nu_n}},1\right[\right)$$

$$= 2B_{n-l-1}^{-1}\int_{1/2+K/\sqrt{nu_n}}^1 x^{n-l-1}(1-x)^{n-l-1}\,dx$$

$$\leq B_{n-l-1}^{-1}\left(\frac{1}{4}-\frac{K^2}{nu_n}\right)^{n-l-1}\left(\frac{1}{2}-\frac{K}{\sqrt{nu_n}}\right).$$

Thus, using (11) and the expansion

$$\ln\left(\frac{1}{4}-\frac{K^2}{nu_n}\right) = -\ln 4 + 1 - \frac{4K^2}{nu_n} + o\left(\frac{1}{nu_n}\right),$$

we obtain

$$\limsup_n u_n \ln \mathbb{P}(X_n \in [-K,K]^c) \leq -4K^2,$$

which implies the exponential tightness. □

LEMMA 4.5. *For every sequence $(u_n)$ decreasing to $0$ with $n^{-1}=o(u_n)$ and $k\in\mathbb{N}^*$, the sequence of random vectors $Z'_n := 2\sqrt{nu_n}A_k(p^{(k)}(Z_n)-\bar{p}^{(k)})$ satisfies a LDP with good rate function $J_k$ and speed $(u_n)$.*

PROOF. Using the two previous lemmas, we can state a LDP for the random vectors $W_n := \sqrt{nu_n}(p^{(k)}(Z_n) - \bar{p}^{(k)}))$ with good rate function

$$J'_k(\bar{x}) = 4\bar{x}^T\bar{x}.$$

So, using the contraction principle, we obtain a LDP for $Z'_n = 2A_k W_n$ with good rate function

$$\begin{aligned}J_k(\bar{x}) &= \inf\{J'_k(\bar{y}) : \bar{x} = 2A_k\bar{y}\} = \tfrac{1}{4}\bar{x}^T(A_k^{-1})^T A_k^{-1}\bar{x} \\ &= \tfrac{1}{2}\bar{x}^T\Sigma_k^{-1}\bar{x}.\end{aligned}\qquad\square$$

LEMMA 4.6. *The random vector sequences $(\widetilde{Z}_n^{(k)})$ and $(Z'_n)$ are exponentially equivalent.*

PROOF. It is obvious that there exists a constant $C_0$ such that

$$|\mathbf{O}(|p^{(k)} - \bar{p}^{(k)}|^2)| \leq C_0 |p^{(k)} - \bar{p}^{(k)}|^2,$$

where $\mathbf{O}(|p^{(k)} - \bar{p}^{(k)}|^2)$ is the function which appears in (9). So,

$$\begin{aligned}\{|\widetilde{Z}_n^{(k)} - Z'_n| > \varepsilon\} &= \{\sqrt{nu_n}|\mathbf{O}(|p^{(k)} - \bar{p}^{(k)}|^2)| > \varepsilon\} \\ &\subset \left\{\sqrt{nu_n}|p^{(k)} - \bar{p}^{(k)}|^2 > \frac{\varepsilon}{C_0}\right\}.\end{aligned}$$

Let $C_1 > 0$. For $\sqrt{nu_n} > C_1$,

$$\left\{\sqrt{nu_n}|p^{(k)} - \bar{p}^{(k)}|^2 > \frac{\varepsilon}{kC_0}\right\} \subset \left\{\sqrt{nu_n}|p^{(k)} - \bar{p}^{(k)}| > \left(\frac{\varepsilon C_1}{kC_0}\right)^{1/2}\right\}.$$

This leads to

$$\limsup_n u_n \ln \mathbb{P}(|\widetilde{Z}_n^{(k)} - Z'_n| > \varepsilon)$$

$$\leq \limsup_n u_n \ln \mathbb{P}\left(\sqrt{nu_n}|p^{(k)} - \bar{p}^{(k)}| > \left(\frac{\varepsilon C_1}{kC_0}\right)^{1/2}\right) \leq -4\frac{\varepsilon C_1}{kC_0},$$

where the last inequality follows from the LDP for the random vectors $W_n$ (Lemma 4.5). Therefore, taking $C_1 \to \infty$, we obtain

$$\limsup_n u_n \ln \mathbb{P}(|\widetilde{Z}_n^{(k)} - Z'_n| > \varepsilon) = -\infty.$$

Then $(\widetilde{Z}_n^{(k)})$ and $(Z'_n)$ are exponentially equivalent. $\square$

Theorem 2.7 is obtained directly using the two previous lemmas.

4.4. *Proof of Theorem 2.3.*

4.4.1. *Statement of the LDP.* Consider on $M$ the product topology and the product algebra. Let $M_n^+$ be the subset of $M$ defined by

$$M_n^+ = \{c^{(\infty)} \in M : c_{n+1} = c^+(c^{(n)})\}.$$

We define on $M$ the following sequence of p.m.s $(\bar{P}_n)_{n \in \mathbb{N}}$:

$$\bar{P}_n(B) = P_n(\Pi_n(B \cap M_n^+)), \qquad B \text{ measurable set of } M.$$

THEOREM 4.7. $(\bar{P}_n)_{n \in \mathbb{N}}$ *satisfies a LDP with rate function*

$$I_M(c^{(\infty)}) = \lim_n -\ln(4^n(c^+(c^{(n)}) - c^-(c^{(n)}))).$$

PROOF. Let $k \in \mathbb{N}$ and $B_k$ be a measurable set of $M_k$. We have, for $n \geq k$,

$$\bar{P}_n \circ \Pi_k^{-1}(B_k) = \bar{P}_n(B_k \times \mathbb{R}^\mathbb{N}) = P_n(\Pi_n(B_k \times \mathbb{R}^\mathbb{N} \cap M_n^+)) = P_n(B_k \times \mathbb{R}^{n-k})$$
$$= \mathbb{P}(Z_n^{(k)} \in B_k).$$

Therefore, the family $\{\bar{P}_n \circ \Pi_k^{-1}\}_{n \in \mathbb{N}}$ satisfies a LDP with rate function $I_k$. The Dawson–Gärtner theorem [Theorem 4.6.1 of Dembo and Zeitouni (1998)] leads to a LDP for $\bar{P}_n$ on the projective limit of the spaces $M_n$ which is $M$. The good rate function is $I_M(c^{(\infty)}) = \sup_n I_n(\Pi_n c^{(\infty)})$. For $c^{(\infty)} \in M$, $I_n(\Pi_n c^{(\infty)})$ is nonincreasing in $n$, thus, $\sup_n I_n = \lim_n I_n$ [see (10)]. □

Let $G$ denote the one-to-one correspondence between $M$ and $\mathbb{P}([0,1])$. It is easy to see that it is continuous. Therefore, using the contraction principle, we obtain a LDP for the probability measure family $\{\bar{P}_n \circ G^{-1}\}_{n \in \mathbb{N}}$. But, for $B$ Borelian of $\mathbb{P}([0,1])$,

$$Q_n(B) = (\bar{P}_n \circ G^{-1})(B) = \bar{P}_n(G^{-1}(B)).$$

We have $c^{(\infty)} \in G^{-1}(B) \cap M_n^+$ if, and only if, $G(c^{(\infty)}) = \sigma_n^+(c^{(n)}) \in B$. So,

$$(\bar{P}_n \circ G^{-1})(B) = \mathbb{P}(Z_n \in \Pi_n(G^{-1}(B) \cap M_n^+))$$
$$= \mathbb{P}(\sigma_n^+(Z_n) \in B).$$

Therefore, the contraction principle gives the LDP for $(Q_n)_{n \in \mathbb{N}}$ with good rate function

(14) $$\tilde{I}(\mu) = \lim_n -\ln[4^{n-1}(c_n^+(\mu) - c_n^-(\mu))],$$

where $c_n^+(\mu)$ and $c_n^-(\mu)$ denote $c^+(c^{(n-1)}(\mu))$ and $c^-(c^{(n-1)}(\mu))$, respectively. Obviously, from (10), we may conclude that $\tilde{I}$ achieves the value 0 only at $\nu$.

LEMMA 4.8. *The function $\tilde{I}$ in* (14) *is convex.*

PROOF. Let $\mu = \alpha\mu_1 + \beta\mu_2$ with $\mu,\mu_1,\mu_2 \in \mathbb{P}([0,1])$ and $\alpha,\beta$ positive real numbers such that $\alpha + \beta = 1$. Thus,
$$c^{(n)}(\mu) = \alpha c^{(n)}(\mu_1) + \beta c^{(n)}(\mu_2).$$
Thanks to convexity of $M_{n+1}$,
$$\alpha c^+(c^{(n)}(\mu_1)) + \beta c^+(c^{(n)}(\mu_2)) \in \{c : (c^{(n)}(\mu), c) \in M_{n+1}\}$$
and by the definition of $c^+$, we have
$$c^+(c^{(n)}(\mu)) \geq \alpha c^+(c^{(n)}(\mu_1)) + \beta c^+(c^{(n)}(\mu_2)).$$
Analogously, we obtain that $c^-(c^{(n)}(\mu)) \leq \alpha c^-(c^{(n)}(\mu_1)) + \beta c^-(c^{(n)}(\mu_2))$. Thus,
$$c_n^+(\mu) - c_n^-(\mu) \geq \alpha(c_n^+(\mu_1) - c_n^-(\mu_1)) + \beta(c_n^+(\mu_2) - c_n^-(\mu_2)).$$
Since $t \mapsto -\ln 4^n t$ is a nonincreasing convex function in $t > 0$ for all $n \in \mathbb{N}^*$, we obtain that $\tilde{I}$ is convex as the supremum of convex functions. □

4.4.2. *Identification of the rate function.* From Corollary 1.4.6 of Dette and Studden (1997), we may write the following equation:

$$(15) \quad r_n(c^{(n-1)}) = \frac{\underline{\mathrm{H}}_{n-1}\bar{\mathrm{H}}_{n-1}}{\underline{\mathrm{H}}_{n-2}\bar{\mathrm{H}}_{n-2}}, \qquad c^{(n-1)} = (c_1,\ldots,c_{n-1})^T \in M_{n-1},$$

where, for $c^{(2m+1)} \in M_{2m+1}$ $(m \in \mathbb{N}^*)$,

$$\underline{\mathrm{H}}_{2m} = \begin{vmatrix} c_0 & \cdots & c_m \\ \vdots & & \vdots \\ c_m & \cdots & c_{2m} \end{vmatrix},$$

$$\bar{\mathrm{H}}_{2m} = \begin{vmatrix} c_1 - c_2 & \cdots & c_m - c_{m+1} \\ \vdots & & \vdots \\ c_m - c_{m+1} & \cdots & c_{2m-1} - c_{2m} \end{vmatrix},$$

$$\underline{\mathrm{H}}_{2m+1} = \begin{vmatrix} c_1 & \cdots & c_{m+1} \\ \vdots & & \vdots \\ c_{m+1} & \cdots & c_{2m+1} \end{vmatrix},$$

$$\bar{\mathrm{H}}_{2m+1} = \begin{vmatrix} c_0 - c_1 & \cdots & c_m - c_{m+1} \\ \vdots & & \vdots \\ c_m - c_{m+1} & \cdots & c_{2m} - c_{2m+1} \end{vmatrix}.$$

These determinants are called *Hankel determinants*. Within the moment problem they play an important role. The conditions under which a sequence is a moment sequence can be expressed on these Hankel determinants [see Karlin and Studden (1966), Kreĭn and Nudelman (1977) and Dette and Studden (1997)].

We have the following relations between Hankel determinants. Let $\mu \in \mathbb{P}([0,1])$. Let $\mu' \in \mathbb{P}([0,1])$ be defined by

$$\mu'(dx) = \frac{x - x^2}{c_1(\mu) - c_2(\mu)} \mu(dx).$$

Obviously, $\mu'$ verifies

$$c_k(\mu') = \frac{c_{k+1}(\mu) - c_{k+2}(\mu)}{c_1(\mu) - c_2(\mu)}.$$

Hence,

$$\bar{H}_{2m}(\mu) = (c_1(\mu) - c_2(\mu))^m \underline{H}_{2(m-1)}(\mu').$$

Moreover,

$$\frac{\bar{H}_{2m}(\mu)}{\bar{H}_{2(m-1)}(\mu)} = (c_1(\mu) - c_2(\mu)) \frac{\underline{H}_{2(m-1)}(\mu')}{\underline{H}_{2(m-2)}(\mu')}.$$

In what follows we set $c^{(n)} = c^{(n)}(\mu)$, $c_n^+ = c^+(c^{(n-1)})$, and $c_n^- = c^-(c^{(n-1)})$. The Hankel matrices depending on $\mu'$ are tagged with a prime. We have

$$(c_{n+1}^+ - c_{n+1}^-)(c_n^+ - c_n^-) = \frac{\bar{H}_n \underline{H}_n}{\bar{H}_{n-2} \underline{H}_{n-2}}$$

$$= (c_1 - c_2) \frac{\underline{H}'_{n-2}}{\underline{H}'_{n-4}} \frac{\underline{H}_n}{\underline{H}_{n-2}},$$

where the last equality is only true if $n$ is even.

As the limit in (14) exists, we can calculate it taking $n$ even,

$$\lim_m \ln(4^{2m-1}(c_{2m}^+ - c_{2m}^-))$$

$$= \frac{1}{2} \lim_m \ln(4^{2m}(c_{2m+1}^+ - c_{2m+1}^-))$$

$$+ \frac{1}{2} \lim_m \ln(4^{2m-1}(c_{2m}^+ - c_{2m}^-))$$

$$= \frac{1}{2} \lim_m \ln(4^{4m-1}(c_{2m+1}^+ - c_{2m+1}^-)(c_{2m}^+ - c_{2m}^-))$$

$$= \frac{1}{2} \lim_m \ln\left(4^{4m-1}(c_1 - c_2) \frac{\underline{H}'_{2(m-1)}}{\underline{H}'_{2(m-2)}} \frac{\underline{H}_{2m}}{\underline{H}_{2(m-1)}}\right).$$

Grenander and Szegö (1958) have showed a general limit theorem for the quotient of Hankel determinants (see Theorem in Section 5.2 and Section 6.3 of this book). In our framework it can be written in the following way.

THEOREM 4.9. *Let $\mu \in \mathbb{P}([0,1])$. Call $f(x)$ its Radon–Nikodym derivative with respect to the Lebesgue measure, then*

$$\lim_n 4^{2n+1} \frac{\underline{\mathrm{H}}_{2n}(\mu)}{\underline{\mathrm{H}}_{2n-2}(\mu)} = 2\pi \exp\bigg(\int_0^1 \ln f(x) \, d\nu(x)\bigg).$$

Now we may write

$$\begin{aligned}
\lim_m \ln(4^{2m-1}&(c_{2m}^+ - c_{2m}^-)) \\
&= \frac{1}{2} \lim_m \ln\bigg(4^{2m+1} \frac{\underline{\mathrm{H}}_{2m}}{\underline{\mathrm{H}}_{2(m-1)}}\bigg) + \frac{1}{2} \lim_m \ln\bigg(4^{2m-1} \frac{\underline{\mathrm{H}}'_{2(m-1)}}{\underline{\mathrm{H}}'_{2(m-2)}}\bigg) \\
&\quad + \frac{1}{2} \lim_m \ln(4^{-1}(c_1 - c_2)) \\
&= \frac{1}{2} \int_0^1 \ln f(x) \, d\nu(x) + \frac{1}{2} \int_0^1 \ln \frac{(x-x^2)f(x)}{c_1 - c_2} \, d\nu(x) \\
&\quad + \frac{1}{2} \ln(c_1 - c_2) + \frac{1}{2} \ln(\pi) \\
&= \frac{1}{2} \int_0^1 \ln[\pi(x - x^2) f^2(x)] \, d\nu(x) = -I(\mu).
\end{aligned}$$

□

4.5. *Proof of Theorem* 2.5. The LDP for $(\widetilde{Q}_n)$ is a direct consequence of Theorem III.17 of den Hollander (2000). The rate function which controls the LDP is

$$I_F(\mu) = I(\mu) - F(\mu) + K_F,$$

where $K_F := \sup_{\mu' \in \mathbb{P}([0,1])} \{F(\mu') - I(\mu')\}$. □

**Acknowledgments.** The second author is very grateful to Professor Didier Dacuhna-Castelle for his help allowing him to spend some months in French Universities. He also thanks Le Laboratoire de Statistique et Probabilités of Toulouse University for its warm hospitality. The authors also thank an anonymous referee, for a careful reading enhanced the first version of this work, and Professor Alain Rouault who has suggested some simplification in the proof of Lemma 4.3.

Laboratoire de Statistique
et Probabilités
Université Paul Sabatier
UMR C5583
118 route de Narbonne
31062 Toulouse cedex 4
France
e-mail: gamboa@cict.fr
url: www.lsp.ups-tlse.fr/Fp/Gamboa/index.html

Laboratoire de Statistique
et Probabilités
Université Paul Sabatier
France
and
Facultad de Matemáticas y Computación
Universidad de la Habana
San Lázaro y L
CP 10400 Ciudad de la Habana
Cuba
e-mail: livang@matcom.uh.cu